\documentclass[12pt]{article}
\usepackage{amssymb}

\def\appendix{\par}  % to see appendix

\def\al{\alpha}
\def\anal{{\bf \Sigma}^1_1}
\def\base{{\mathcal B}}
\def\bb{{\mathfrak b}}
\def\be{\beta}
\def\bool{{\mathbb B}}
\def\borel{{\rm Borel}}
\def\bsn{{\mathfrak s \mathfrak n}^*}
\def\coanal{{\bf \Pi}^1_1}
\def\de{\delta}
\def\diff{{\Delta}}
\def\ff{{\mathcal F}}

\def\gak{{\gamma_k}}
\def\ga{\gamma}
\def\gg{{\mathcal G}}
\def\ka{\kappa}
\def\nonmeag{{\rm non}({\mathcal M})}
\def\om{\omega}
\def\ord{{\rm ord}}
\def\proof{\par\noindent Proof\par\noindent}
\def\pr{\prime}
\def\qed{\par\noindent QED\par}
\def\rationals{{\mathbb Q}}
\def\reals{{\mathbb R}}
\def\res{\upharpoonright}
\def\si{\sigma}
\def\sm{{\setminus}}
\def\sn{{\mathfrak s \mathfrak n}}
\def\sq{\subseteq}
\def\st{\;:\;} % such that
\def\uu{{\mathcal U}}

\newtheorem{theorem}{Theorem}[section]
\newtheorem{lemma}[theorem]{Lemma}

\newtheorem{cor}[theorem]{Corollary}

\newtheorem{prop}[theorem]{Proposition}
\newtheorem{example}[theorem]{Example}

\begin{document}

\begin{center}
{\large A hodgepodge of sets of reals }
\end{center}

\begin{flushright}
Arnold W. Miller\footnote{
Thanks to the conference organizers: Cosimo Guido, Ljubisa Ko\v{c}inac, Boaz
Tsaban, Liljana Babinkostova, and Marion Scheepers for their generosity in
inviting me to speak at the Workshop on Coverings, Selections and Games in
Topology, December 2005, University of Lecce, Italy. 
 }
\end{flushright}

\begin{center}
Abstract
\end{center}

\begin{quote}
We open up a grab bag of miscellaneous results and remarks about sets
of reals.  The reader
will be treated to a cornucopia of delightful, delectable, and delicious
ideas
to pounder and shake his head at in consternation while muttering
``Who would have thought of that?'' and ``Why didn't they keep it
to themselves?''.  Results concern: Kysiak and Laver-null sets,
Ko\v{c}inac
and $\gak$-sets, Fleissner and square $Q$-sets,  Alikhani-Koopaei and
minimal $Q$-like-sets, Rubin and $\sigma$-sets, and
Zapletal and the Souslin number.
See the survey papers Brown, Cox \cite{browncox},
and Miller \cite{survey1,survey2}.
\end{quote}

\section{$\si$-sets are Laver null}

A subtree $T\sq\om^{<\om}$ of the finite sequences
of elements of $\om=\{0,1,2,\ldots\}$ is called a Laver tree
(\cite{laver}) iff there exists $s\in T$ (called the root node of
$T$) with the property that for every $t\in T$ with $s\sq t$ there
are infinitely many $n\in\om$ with $t n$ in $T$.  Here $tn$ is
the sequence of length exactly one more than $t$ and ending in
$n$.  We use $[T]$ to denote the infinite branches of $T$, i.e.,
$$[T]=\{x\in\om^\om\st \forall n\in\om\;\; x\res n\in T\}$$

A set $X\sq \om^\om$ is Laver-null iff for every Laver tree
$T$ there exists a Laver subtree $T^\pr\sq T$ such that
$$[T^\pr]\cap X=\emptyset$$
This is analogous to the ideal of Marczewski null sets, $(s)_0$.
For some background on this topic,
see Kysiak and Weiss \cite{kysiak} and Brown \cite{brown}.

A separable metric space $X$ is a $\si$-set iff every
$G_\de$ in $X$ is also $F_\si$.  It is known to be
relatively consistent
(Miller \cite{borhier})
with the usual axioms of set theory that every $\si$-set
is countable.

At the Lecce conference, Kysiak\footnote{
As I was writing this I learned from Jack Brown that M.Kysiak, A.Nowik, and
T.Weiss \cite{KNW},
also solved this problem at about the same time.  In fact, their
solution is a little better as it also solves the analogous problem for Ramsey
null sets.}
asked if it is consistent to have a $\si$-set
which is not Laver-null. The answer is no.

\begin{theorem}
Every $\si$-set is Laver-null.  In fact, the
Borel hierarchy of a non-Laver-null set must have $\om_1$ levels.
\end{theorem}

\proof
Here we use a result of Rec\l aw that appears in Miller \cite{survey2}.
Rec\l aw proved that if $X$ is a set of reals and there exists
a continuous onto map
$f:X\to 2^\om$, then the Borel hierarchy on $X$
has $\om_1$ levels, in particular, $X$ is not a $\si$-set.

So let $X\sq\om^\om$ be a set which is not Laver-null.  Hence
there exists a Laver tree $T$ such that for every Laver subtree
$T^\pr\sq T$ we have that $[T^\pr]$ meets $X$.

To simplify our
notation assume that $T=\om^{<\om}$.  Define the following continuous
function $f:\om^\om\to 2^\om$.  $f$ is the parity function, i.e.,
for any $x\in \om^\om$ we define $f(x)=y\in 2^\om$ by the
rule that $x(n)$ is even iff $y(n)=0$.  But note that
$f$ maps $X$ continuously onto $2^\om$.  This is because for
any $y\in 2^\om$ there is a Laver-tree $T$ such that
$f[T]=\{y\}$.  But since $[T]$ meets $X$ there is some $x\in X$
with $f(x)=y$.

Now in the more general case $T$ is an arbitrary Laver-tree. In
this case note that there is a natural map from $\om^{<\om}$ to
$T$ and by using essentially the same proof as above it is easy
to see that the result holds.

\qed

\section{$\gak$-sets}

In  Ko\v{c}inak \cite{koc} the notion of a $\gak$-set is defined.
A $k$-cover of topological space $X$ is a family of
open subsets with the property that every compact subset of
$X$ is subset of an element of the family.  $X$ is called
$\gak$-set iff for every $k$-cover $\uu$ of
$X$ there exists a sequence $(U_n\in\uu : n\in\om)$ such that
for every compact $C\sq X$ we have that $C\sq U_n$ for all but finitely
many $n$.

This is a generalization of $\ga$-sets first considered by
Gerlits-Nagy \cite{gn} and studied in many papers.

A theorem of Galvin and Todor\v cevi\'c 
(see Galvin and Miller \cite{gm}) shows that
it is consistent that the union of two $\ga$-sets need not be a $\ga$-set.
Ko\v{c}inak asked at the Lecce conference if such a counterexample exists for
$\gak$-sets.  We show that it does.

\begin{example}\label{exampkoc}
There exist disjoint subsets of the plane $X$ and $Y$ such
that both $X$ and $Y$ are $\gak$-sets but
$X\cup Y$ is not.
\end{example}

Let $X$ be the open disk of radius one, i.e.,
$X=\{(x,y)\st x^2+y^2<1\}$,
and $Y$ be any singleton on the boundary of $X$, e.g.,
$Y=\{(1,0)\}$.
The result follows easily from the following:

\begin{lemma}
Suppose that $Z$ is a metric space.  Then $Z$ is
a $\gak$-set iff $Z$ is locally compact and separable.
\end{lemma}
\proof
First suppose that $Z$ is locally compact and separable.  Then
we can write $Z$ as an increasing union of compact subsets
$C_n$ whose interiors cover $Z$.  
Given a $k$-cover $\uu$ we simply choose $U_n\in\uu$
so that $C_n\sq U_n$.  This works because for every compact set $C$
there exists $n$ with $C\sq C_n$.

Conversely, suppose that $Z$ is not locally compact. This
means that for some $x\in Z$ we have that $x$ is not in the
interior of any compact set.  Define a sequence of
$\uu_n$ as follows:
Let $\uu_n$ be the set of all open subsets of $Z$ such
that $U$ does not contain the open ball of radius $1/2^n$ around
$x$, i.e. there exists $y\notin U$ such that $d(x,y)<1/2^n$.

Note that each $\uu_n$ is a $k$-cover of $Z$.  To see this, suppose $C$
is a compact subset of $Z$.  Since $x$ is not in the interior of
$C$, the set $C$ cannot contain an open ball centered at $x$.
Choose $y\notin C$ with $d(x,y)<1/2^n$.  Now cover $C$ with
(finitely) many open balls not containing $y$. The union of this
cover is in $\uu_n$.

We can use the trick of Gerlits and Nagy to get a single $k$-cover
from the
sequence  of
$k$-covers, $(\uu_n\st n\in\om)$.
Since $Z$ cannot be compact there
must exist a sequence $(x_n:n\in\om)$ with no limit point.  Define
$$\uu=\{U\sm\{x_n\} \st n<\om,\;\; U\in \uu_n\}.$$
Since any compact set can contain at most finitely many of
the $x_n$, we see that $\uu$ is a $k$-cover of $Z$.

For contradiction, suppose $Z$ is $\gak$-set and
$(U_n\in\uu:n\in\om)$ eventually contains each compact set.
Without loss, we may assume that $U_n\in \uu_{l_n}$ with
$l_n$ distinct. This is because at most finitely many $U_n$ can
be ``from'' any $\uu_l$ since they eventually must include $x_l$.
Choose $y_n\notin U_n$ with $d(x,y_n)<1/{2^{l_n}}$.  Then
$$\{y_n\st n\in\om\}\cup \{x\}$$
is a convergent sequence, hence compact. But it is not a
subset of any $U_n$.

It is easy to see that $Z$ must be separable as we can take
$\uu_n$ to be the family of finite unions of open balls of
radius less than $1/2^n$, then apply the Gerlits Nagy trick
as above to obtain a countable basis for $Z$.

\qed

In Example \ref{exampkoc}
each of $X$ and $Y$ are locally compact metric spaces but
$X\cup Y$ is not locally compact at the point $(0,1)$, so the
result follows.

Ko\v{c}niac also asked if $X\times Y$ is $\gak$-set if both
$X$ and $Y$ are.  For metric spaces, this must be true by the Lemma,
since the product of locally compact separable metric spaces is a locally
compact separable metric space.

\section{Q-sets}

A Q-set is a separable metric space $X$ such that every subset
of $X$ is a (relative) $G_\de$-set.  It is easy to
see that $2^{|X|}=2^\om$, hence, if there is an uncountable
Q-set, then $2^{\aleph_1}=2^{\aleph_0}$.  So uncountable Q-sets
might not exist.  Martin's axiom (MA) implies that every separable
metric space of size less than the continuum is a Q-set
(see Martin and Solovay \cite{internal}).

The Rothberger cardinal, $\bb$, is defined to be the cardinality
of the smallest family $\ff\sq\om^\om$ such that for every $g\in\om^\om$
there is some $f\in \ff$ with $f(n)\geq g(n)$ for infinitely many $n$.
That is to say, $\bb$ is the size of the smallest unbounded family
in the quasi-ordering $(\om^\om,\leq^*)$.
Martin's Axiom implies that $\bb$ is the continuum.

\begin{theorem} \label{qsets}
Suppose $\ka<\bb$.  Then the following are equivalent:
\begin{enumerate}
\item \label{existq} There exists a Q-set $X\sq 2^\om$ with $|X|=\ka$.
\item \label{functionq}
There exists $(f_\al:\om^\om\to 2^\om\st \al<\ka)$ continuous
functions such that given any $(y_\al\in 2^\om\st \al<\ka)$ there exists
$x\in \om^\om$ with the property that
$f_\al(x)=^*y_\al$ for every $\al<\ka$.
\item \label{univq}
There exists a sequence $(U_\al\sq 2^\om\times 2^\om\st\al<\ka)$
of $G_\de$-sets which is universal for $\ka$ sequences
of $G_\de$-sets, i.e., for every sequence
$$(V_\al\sq 2^\om\st\al<\ka)$$
of $G_\de$-sets there exists $x\in 2^\om$ such that for every $\al<\ka$
$$V_\al=U_\al(x)=^{def}\{y\st (x,y)\in U_\al\}.$$
\end{enumerate}
\end{theorem}
\proof

We will need the following lemma and the details of its proof.

\begin{lemma} \label{lemq}
There exists $U\sq 2^\om\times 2^\om$ which is a universal $G_\de$-set
such that for every $x_1,x_2\in 2^\om$ if
$x_1=^*x_2$, then $U(x_1)=U(x_2)$.
\end{lemma}
\proof
A set $U$ is a universal $G_\de$-set, if it is $G_\de$ and for every
$G_\de$-set $V\sq 2^\om$ there exists $x\in 2^\om$ such that
$$U(x)=^{def}\{y\in 2^\om\st (x,y)\in U\}=V.$$
Define
$$U=\{(A,y)\in P(2^{<\om})\times 2^\om\st\exists^\infty n\;\; y\res n\in A\}$$
where $\exists^\infty$ stands for ``there exists infinitely many''.
It is easy to see that $U$ is $G_\de$.   To see that it is universal, suppose
that $V=\cap_{n<\om} V_n$ where the $V_n\sq 2^\om$ are open and
descending, i.e.,  $V_{n+1}\sq V_n$ for each $n$.
For $\si\in 2^{<\om}$ nontrivial
let $\si^*\sq\si$ be the initial
segment of $\si$ of length exactly one less than $\si$,
i.e., $|\si^*|=|\si|-1$.   Define
$$A=\{\si\st [\si]\sq V \mbox{ or }\exists n\;\; [\si]\sq V_n
\mbox{ and } [\si^*]\not\sq V_n\}$$
Then $U(A)=V$.  To see this,  suppose $x\in U(A)$.
If for some $n$ we have that $x\res n\in A$ because $[x\res n]\sq V$
then clearly $x\in V$.  On the other hand, if there are infinitely
many $k$ such that for some $n$, $[x\res k]\sq V_n$ but
$[x\res (k-1)]\not\sq V_n$, then these $n$'s must all be distinct
and since the $V_n$ were descending $x\in V$.

Conversely, if $x\in V$
then either $x$ is in the interior of $V$ and so
$x\res k\in A$ for all but finitely many $k$ or it isn't in the
interior of $V$ and there are thus infinitely many $n$ with
$x\res n\in A$.  Hence $x\in U(A)$.

\qed

\bigskip
\ref{functionq}$\to$\ref{univq}:
This follows immediately from the Lemma.  Just define
$$(x,y)\in U_\al\mbox{ iff } (f_\al(x),y)\in U$$
identify $\om^\om$ with a $G_\de$ subset of $2^\om$.

\bigskip
\ref{univq}$\to$\ref{existq}:
By the proof of Lemma \ref{lemq} there exists
$A_\al\sq 2^{<\om}\times 2^{<\om}$ such that
for any $(x,y)$ we have that
$(x,y)\in U_\al$ iff $\exists^\infty n\;\; (x\res n,y\res n)\in A_\al$.
We claim that
$$\{A_\al\st \al<\ka\}$$
is a $Q$-set.  Fix $y\in 2^\om$ arbitrary.  Consider any
$\Gamma\sq \ka$ and define the sequence of $G_\de$ sets
$(V_\al:\al<\ka)$ by 
$$V_\al=\left\{
\begin{array}{ll}
\{y\} & \mbox{ if } \al\in\Gamma\\
\emptyset & \mbox{ if } \al\notin\Gamma\\
\end{array}\right.
$$ 
By
assumption there exists $x\in 2^\om$ such that
$U(\al)=V_\al$ for every $\al<\ka$.  But then
$$\al\in\Gamma \mbox{ iff }
y\in U_\al(x) \mbox{ iff }
\exists^\infty n\;\; (x\res n,y\res n)\in A_\al \mbox{ iff } $$
$$A_\al\in \{A\st \exists^\infty n\;\;(x\res n,y\res n)\in A\}$$
But this last set is $G_\de$.  It follows
that $\{A_\al\st\al\in\Gamma\}$ is relatively $G_\de$
in the set $\{A_\al\st\al\in\ka\}$.

\bigskip
\ref{existq} $\to$ \ref{functionq}:
Let $\{v_\al^n\in 2^\om \st n<\om,\;\;\al<\ka\}$
be a $Q$-set.  Define the function $f_\al$ as follows:
Suppose $x=(A,(I_n:n<\om))$ where $A\sq 2^{<\om}$ and
each $I_n\sq 2^{<\om}$ is finite.  (We can easily identify the set
of such $x$ with
$\om^\om$.)   Now for each $\al<\ka$ define a continuous map
$f_\al(x)\in 2^\om$ as follows.   Define

$$f_\al((A,(I_n:n<\om)))(n)=
\left\{
\begin{array}{ll}
1 & \mbox{ if } \exists k\;\; v_\al^n\res k\in I_n\cap A\\
0 & \mbox{ otherwise }\\
\end{array}\right.
$$ 
Since the $I_n$ are finite, the function
$f_\al$ is continuous.  We verify that it has the property required.
Let $x_\al\in 2^\om$ for $\al<\ka$ be arbitrary.  Since
$\{v_\al^n\in 2^\om \st n<\om,\;\;\al<\ka\}$ is a $Q$-set,
there is a $G_\de$-set $U\sq 2^\om$ with
the property that for every $\al<\ka$ and $n<\om$ we
have that $v_\al^n\in U$ iff $x_\al(n)=1$.
By the proof of Lemma \ref{lemq} there exists $A\sq 2^{<\om}$ such
that for all $\al,n$

$v_\al^n\in U$ iff $A\cap\{v_\al^n\res k\st k<\om\}$ is
infinite.

\noindent Since $\bb>\ka$ there exists a partition
$(I_l : l<\om)$ of $2^{<\om}$ into finite sets such that
for every $\al<\ka$ and $n<\om$
the set $A\cap\{v_\al^n\res k\st k<\om\}$ is
infinite iff $I_l\cap A\cap\{v_\al^n\res k\st k<\om\}\neq \emptyset$
for all but finitely many $l<\om$.  But this implies that
for $f_\al((A,(I_l : l<\om))=^*x_\al$ for each $\al$.

\qed

\bigskip
Condition \ref{univq} is a kind of uncountable version of Luzin's doubly
universal sets, see Kechris \cite{kec} page 171 22.15 iv.  Luzin used
a doubly universal set to prove that the classical properties
of separation and reduction cannot hold on the same side of a
reasonable point-class.

In condition \ref{functionq},  $u=^*v$
means that $u(n)=v(n)$ except for finitely many $n$.  It is impossible
to have the stronger condition with ``$=$'' in place of ``$=^*$'' at least
when $\ka$ is uncountable.  To see this, fix $y_0\in 2^\om$ and
define $E_\al=f_\al^{-1}(y_0)$ for $\al<\om_1$.  It is not
hard to see that the $F_\al=\cap_{\be<\al}E_\be$ would have
to be a strictly decreasing sequence of closed sets, which is impossible
in a separable metric space.

We do not know iff the condition $\ka<\bb$ is needed for this result.
There are several models of set theory where there is a $Q$-set
and $\bb=\om_1$,  Fleissner and Miller \cite{FM},
Judah and Shelah \cite{JS}, and Miller \cite{madq}.

We obtained this result while working on the square Q-set problem,
see Fleissner \cite{squareQ}.  Unfortunately, Fleissner's proof that
it is consistent there is a $Q$-set whose square is not a $Q$-set contains
a gap. In his paper:
he claims to show that in his model of set theory:
\begin{enumerate}
\item there is a $Q$-set $Y\sq 2^\om$ of size $\om_2$, and
\item for any set of $Z=\{z_\al \st \al< \om_2\}\sq  2^\om$ the set
$$\{(z_\al,z_\be) \st \al<\be<\om_2\}$$
 is not $G_\de$ in $Z \times Z$.
\end{enumerate}
But we have a fairly easy proof that (1) implies the negation of (2).

\begin{theorem}
If there exists a Q-set $Y\sq 2^\om$
with $|Y|=\om_2$ then there exists
$Z=\{z_\al:\al<\om_2\}\sq 2^\om$ such that
$$\{(z_\al,z_\be):\al<\be<\om_2\}$$
is (relatively) $G_\delta$ in $Z^2$.
\end{theorem}

\proof Let $Y=\{y_\al:\al<\om_2\}$ and let
$U\sq 2^\om\times 2^\om$ be a universal $G_\de$-set.
Choose for each $\be<\om_2$ a $u_\be\in 2^\om$ such that
for every $\alpha<\om_2$
$$y_\alpha \in U_{u_\be}\mbox{ iff } \alpha<\be .$$
Since $U$ is $G_\delta$ there are clopen $C_{n,m},D_{n,m}\sq 2^\om$
with
$$U=\bigcap_{n<\om}\bigcup_{m<\om} (C_{n,m}\times D_{n,m})$$

Now let $z_\alpha=(y_\alpha,u_\alpha)$ and identify $2^\om\times 2^\om$
with $2^\om$.

Then for any $\alpha,\be<\om_2$ we have that

$\alpha<\be$

iff  $(y_\alpha,u_\be)\in U$

iff $(y_\alpha,u_\be)
\in \bigcap_{n<\om}\bigcup_{m<\om} (C_{n,m}\times D_{n,m})$

iff $(z_\alpha,z_\be)=((y_\alpha,u_\alpha),(y_\be,u_\be))
\in
\bigcap_{n<\om}\bigcup_{m<\om} ((C_{n,m}\times 2^\om)\times
(2^\om\times D_{n,m}))$

\qed

As far as we know, the problem of the consistency
of a $Q$-set whose square is not a $Q$-set, is open.
One way to connect this problem with Theorem \ref{qsets} is the
following:

\begin{cor}
Suppose there is a $Q$-set of size $\om_2$ and
$\bb>\om_2$.  Then given any family $\Gamma\sq P(\om_2\times\om_2)$
with $|\Gamma|=\om_2$
there is a $Q$-set $$Z=\{z_\al\in 2^\om\st\al<\om_2\}$$
such that for every $A\in\Gamma$ the
set $\{(z_\al,z_\be)\st (\al,\be)\in A\}$ is
$G_\de$ in $Z$.
\end{cor}

\section{Minimal $Q$-like-sets}

At the Slippery-Rock conference in June 2004, Ali A. Alikhani-Koopaei
asked me if the following $Q$-like example was possible.  We show that
it is.

\begin{example}
There exist a $T_0$ space $Y$ such that $Y$ is not a $Q$-set
but for every
$A\sq Y$ there is a minimal $G_\de$ set $Q$ with
$A\sq Q$.
By minimal we mean that for any $G_\de$ set $Q^\pr$ if
$A\sq Q^\pr$, then $Q\sq Q^\pr$.
\end{example}

\proof
Let $X$ be any Q-set, i.e., every subset of $X$ is $G_\de$ and
$X$ at least $T_0$.  For example, a
discrete space.  Now let $X^\pr$ be a disjoint copy of $X$ and
let $p\mapsto p^\pr$ a bijection from $X$ to $X^\pr$.  For each $A\sq X$
let $A^\pr=\{p^\pr:p\in A\}$.  Define the topology on $Y=X\cup X^\pr$
by letting the open sets of $Y$ be exactly those of the form $U\cup V^\pr$
where $U,V\sq X$ are open in $X$ and $U\sq V$.  Then $Y$ is $T_0$, e.g.
$X^\pr$ is open in $Y$ and separates any $p$ and $p^\pr$.

\bigskip

Claim: For $A,B\sq X$ the set $A\cup B^\pr$ is $G_\de$ in $Y$ iff
$A\sq B$. Furthermore, given any $A,B\sq X$ the set
$A\cup (A\cup B)^\pr$ is the minimal $G_\de$ in $Y$ containing
$A\cup B^\pr$.
\proof

Suppose that $A$ is not a subset of $B$ and let $p\in A\sm B$.  Then any
open set in $Y$ which contains $A$ must also contain $p^\pr$.  The same is
true for any $G_\de$ and hence $A\cup B^\pr$ is not $G_\de$.

On the other hand, suppose $A\sq B$.  Let $A=\cap_{n<\om}U_n$ and
$B=\cap_{n<\om}V_n$ where the $U_n$ and $V_n$ are open in $X$. Now since
$A\sq B$ we may assume that $U_n\sq V_n$ (if not just replace $U_n$ by
$U_n\cap V_n$.  But then in
$$A\cup B^\pr=\cap_{n<\om} (U_n\cup V_n^\pr).$$

For the furthermore, note that if $C\cup D^\pr$ is $G_\de$ and
contains $A\cup B^\pr$, then $A\sq C$ and $C\cup B\sq D$ and so
$A\cup (A\cup B)^\pr\sq C\cup D^\pr$

\qed

\bigskip

Question. Can we get an example which is uncountable but contains no
uncountable Q-set?

\bigskip
Yes.  Let $X=\om_1$ have the topology
with $U\sq X$ is open iff $U=\emptyset$ or
there exists $\al$ with
$$U=[\al,\om_1)=^{def}\{\be:\al\leq\be<\om_1\}$$
Given any $A\sq X$ the smallest $G_\de$ containing $A$
is $[\min(A),\om_1)$.

\section{$\si$-sets and retractive boolean algebras}

The definition of thin set of reals is due to Rubin \cite{rubin} who showed it
equivalent to a certain construction yielding a retractive boolean algebra
which is not the subalgebra of any interval algebra.  Rubin asked whether or
not there is always an uncountable thin set of reals.  We show
that every thin set is a $\si$-set and so by the results of
Miller \cite{borhier} that it is consistent there are no
uncountable $\si$-sets, it is also consistent there are no uncountable
thin sets.

A thin set of reals is defined as follows.  An OIT (ordered interval
tree) is a family of $(G_n:n\in \om)$ such that
each $G_n$ is a family of pairwise disjoint open intervals
such that for $n$ and $I\in G_{n+1}$ there
exists $J\in G_n$ with $I\sq J$.
A set of reals $Y$ is
$(G_n:n\in \om)$-small iff there exists $(F_n\in [G_n]^{<\om}:n\in \om)$
such that for every $x\in Y$ and $n\in\om$ if
$x\in \cup G_n$, then $x\in\cup F_n$.
A set of reals $X$ is thin iff for every OIT $(G_n:n\in \om)$
the set $X$ is a countable union of $(G_n:n\in \om)$-small sets.

\begin{prop}
 If $X\sq\reals$ is thin, then $X$ is a $\si$-set.
\end{prop}
\proof

A thin set cannot contain an interval (see Rubin \cite{rubin})
so without loss of
generality we may suppose that $X$ is disjoint from the rationals
$\rationals$.  Let $\base$ be the family of nonempty open intervals with
end points from $\rationals$.
The following claim is easy to prove and left to the reader.

\bigskip

{\bf Claim}. Given any open set $U\sq I$ where $I\in\base$
we can construct a family of
pairwise disjoint intervals $G\sq \base$ so that
\begin{enumerate}
\item $cl(J)\sq I$ for each
$J\in G$ and
\item $\bigcup G \sq U \sq \bigcup G \cup\rationals$.
\end{enumerate}

\bigskip
Now suppose that $\cap_{n<\om} U_n$ is an arbitrary
$G_\de$ set of reals where
the $U_n$ are open sets. Using the claim it is easy to construct a sequence
$G_n\sq \base$ of pairwise disjoint rational intervals such that:
\begin{enumerate}
\item if $I\in G_{n+1}$, then for some $J\in G_n$ we have $cl(I)\sq J$ and
\item $\bigcup G_n \sq U_n \sq \bigcup G_n \cup\rationals$.
\end{enumerate}
Since $X$ is thin, we have
that $X=\cup_{m<\om}X_m$ where each
$X_m$ is $\{G_n:n<\om\}$-small.
Fix $m$.
There exists $F_{n,m}\in [G_n]^{<\om}$ for $n<\om$ which
witness the smallness of $X_m$.
Let
$$C_m=\cap_{n<\om}(\cup F_{n,m})$$
Note that we may assume that for each $n$ and $I\in F_{n+1,m}$ there
is a $J\in F_{n,m}$ with $cl(I)\sq J$.  Hence
$$\cap_{n<\om}(\cup F_{n,m})=\cap_{n<\om}(\cup_{I\in F_{n,m}}cl(I))$$
and since each
$F_{n,m}$ is finite, $C_m$ is closed.
Since $X$ is disjoint from $\rationals$  we have that
$$X\cap (\cap_{n<\om} U_n) = X\cap (\cup_{m<\om} C_m).$$
Since we started with an arbitrary $G_\de$ set we have that
$X$ is a $\si$-set.
\qed

It is not hard to see that a set of reals is thin iff it is hereditarily
Hurewicz.  See Miller and Fremlin \cite{fremlin} for the definition of
the Hurewicz property.

\section{Souslin number and nonmeager sets}

We obtained these results in March 2004.
First we define the following small cardinal number:

$$\nonmeag=\min\{|X|\;:\; X\sq 2^\om\mbox{ nonmeager }\}$$

For $X\sq 2^\om$ we define $\ord(X)$  (the Borel order of $X$) to be the
smallest $\al<\om_1$ such that every Borel subset $A$ of $2^\om$ there
exist a ${\bf \Sigma}^0_\al$ subset $B$ of $2^\om$ such that
$A\cap X=B\cap X$, if there is no such $\al<\om_1$, we define
$\ord(X)=\om_1$.

To prove our main result (Theorem \ref{two})
we will use the following theorem:

\begin{theorem}\label{one}
There exists $X\sq 2^\om$ with $|X|\leq\nonmeag$ and $\ord(X)=\om_1$.
\end{theorem}
\proof
This is similar to the proof of Miller \cite{borhier} Theorem 18.  Notice
that it is enough to show that for each $\al<\om_1$ there
exists an $X_\al\sq 2^\om$ with
$$|X_\al|\leq \nonmeag$$
 and
$\ord(X_\al)\geq\al$, since the $\om_1$ union of these sets would be
the $X$ we need.

So fix $\al_0<\om_1$ with $\al_0>1$.
According to Miller \cite{borhier} Theorem 13, there exists
a countable subalgebra $\gg\sq \bool$ where $\bool$ is the complete boolean
algebra:
 $$\bool={{\mbox{Borel}(2^\om)}\over {\mbox{meager}(2^\om)}}$$
such that $\gg$ countably generates  $\bool$ in exactly $\al_0$
steps. This last statement means the following:

Define $\gg_0=\gg$. For $\al>0$ an even ordinal define $\gg_\al$ to be the
family of countable disjuncts of elements from $\bigcup_{\be<\al}\gg_\be$ and
for $\al$ an odd ordinal define $\gg_\al$ to be the family of countable
conjuncts of elements from  $\bigcup_{\be<\al}\gg_\be$.   These classes are
analogous to the ${\bf\Sigma}^0_\al$ and ${\bf\Pi}^0_\al$ families of Borel
sets.  Then $\gg$ has the property that $\gg_{\al_0}=\bool$ but for each
$\be<\al_0$, $\gg_\be\not=\bool$.

Now let $Y\sq 2^\om$ be such that $|Y|=\nonmeag$ and  $Y\cap U$ is
nonmeager for every nonempty open subset $U$ of $2^\om$.  Note that
$Y$ has the property that for any Borel subsets $A$ and $B$ of
$2^\om$, if $A\cap Y=B\cap Y$, then the symmetric difference, 
$A\diff B$ is meager.

Let $\ff\sq \borel(2^\om)$
be a family of representatives for $\gg$, i.e.,
 $$\gg=\{[A]:A\in \ff\}$$
where $[A]\in\bool$ is the equivalence class of $A$ modulo the meager 
ideal
in $2^\om$.  Assume $\ff$ is chosen so that the map $A\mapsto [A]$ is
one-to-one and $2^\om$ and $\emptyset$ are the representatives of
$1$ and $0$.
By throwing out a meager subset
of $Y$ we may assume that for any $A,B,C\in \ff$
\begin{enumerate}
\item $[A]\vee [B]=[C]$ iff $(A\cap B)\cap Y=C\cap Y$
\item $[A] \wedge [B]=[C]$ iff $(A\cup B)\cap Y=C\cap Y$
\end{enumerate}
Define $\ff^Y=\{Y\cap A: A\in\ff\}$.
Then we have that $\gg$ and $\ff^Y$ are isomorphic as boolean algebras:
$$(\gg,\vee,\wedge, 0,1)\simeq (\ff^Y,\cup,\cap,\emptyset,Y)$$

Define
$\ff_\be$ and $\ff_\be^Y$ exactly as we did $\gg_\be$ but using countable
unions and intersections instead of disjuncts and conjuncts as we do in a
boolean algebra.

\bigskip
\noindent {\bf Claim}.
\par
\begin{enumerate}
\item By induction on $\be$
\begin{enumerate}
\item $\gg_\be=\{[B]:B\in \ff_\be\}$ and
\item $\ff^Y_\be=\{B\cap Y:B\in \ff_\be\}$.
\end{enumerate}
\item If $\be<\al_0$, then $\ff^Y_\be\not=\ff^Y_{\al_0}$.
\end{enumerate}
\proof
Item (1) is an easy induction.  To see (2)
suppose that $[B]\in\gg_{\al_0}\sm \cup_{\be<\al_0}\gg_{\be}$. 
Without loss $B\in\ff_{\al_0}$ and we claim that
$B\cap Y\in \ff_{\al_0}^Y\sm \cup_{\be<\al_0}\ff_{\be}^Y$. 
Suppose for contradiction that
$B\cap Y\in \ff^Y_\be$ for some $\be<\al_0$.  Then there would exist
$C\in\ff_\be$ with $B\cap Y=C\cap Y$.  But this would imply that
$[B]=[C]\in \gg_\be$ which is a contradiction.
This proves the claim.

\bigskip Now let $\ff^Y=\{C_n:n<\om\}$ and let $i:Y\to 2^\om$ be the
Marczewski characteristic function of the sequence, which is defined
by
 $$i(a)(n)=\left\{
\begin{array}{ll}
1 & \mbox{ if } a\in C_n \\
0 & \mbox{ if } a\notin C_n\\
\end{array}\right.
$$ 
Let $X=i(Y)$.  $i$ need not be one-to-one but by definition, it is onto
$X$, 
so $|X|\leq |Y|=\nonmeag$.  Note that
 $$\{C\cap X\;:\; C \mbox{ is clopen in } 2^\om\}=\{i(C):C\in\ff^Y\}$$
Hence, since the Borel order of $\ff^Y$ is at least $\al_0$ we have that
$\ord(X)\geq\al_0$. This proves Theorem \ref{one}.

\qed

We define the Souslin number $\sn$:
$$\sn=\min\{|X|\;:\; X\sq 2^\om,\;\; 
\exists A\in\anal\;\forall B\in\coanal
\; A\cap X\not= B\cap X\}$$

In Zapletal \cite{zap2} Appendix C, it is shown that
$\sn\geq \bb$, where $\bb$ is the smallest cardinality of
an unbounded family in $\om^\om$.
In Miller \cite{relanal} it is
shown to be consistent to have $\sn > \bb$.

Define the following variant of the Souslin number $\bsn$:
$$\bsn=\min\{|X|\;:\; X\sq 2^\om,\;\; 
\exists A\in\anal\;\forall B\in\borel
\; A\cap X\not= B\cap X\}$$

The following theorem partially confirms a conjecture of Zapletal
that $\sn\leq\nonmeag$, since $\bsn\leq\sn$.   Zapletal was motivated by
results in \cite{zap2} Appendix C and \cite{zap},
which roughly speaking show that it is
impossible to force  $\sn>\nonmeag$ using a 
countable support iteration of
definable real forcing in the presence of suitable large cardinal axioms.
Zapletal's conjecture remains open.

\begin{theorem}\label{two}
$\bsn\leq\nonmeag$
\end{theorem}
\proof

Let $U\sq 2^\om\times 2^\om$ be a universal ${\bf \Sigma}^1_1$ set and 
consider the set of reals $X$
from Theorem \ref{one}.  For each $\al<\om_1$ let $B_\al\sq 2^\om$ be a
${\bf\Sigma}^0_\al$ such that for every $C$ which is ${\bf \Pi}^0_\al$
we have that
 $$B_\al\cap X\not=C\cap X$$
Since $U$ is universal there exists $a_\al\in 2^\om$ such that the cross
section $U_{a_\al}=B_\al$.  Let $Z$ be defined by
 $$Z=\{a_\al:\al<\om_1\}\times X\sq 2^\om\times 2^\om$$
Then $|Z|\leq\nonmeag$ and there is no Borel set $B\sq 2^\om\times 2^\om$ such
that $Z\cap U=Z\cap B$.  This is because if $B$ is say ${\bf \Pi}^0_\al$, then
every cross section of $B$ is ${\bf \Pi}^0_\al$, but then
  $$B_\al\cap X = U_{a_\al}\cap X=B_{a_\al}\cap X$$
which contradicts our choice of $B_\al$.
\qed

\begin{flushleft}
Arnold W. Miller \\
miller@math.wisc.edu \\
http://www.math.wisc.edu/$\sim$miller\\
University of Wisconsin-Madison \\
Department of Mathematics, Van Vleck Hall \\
480 Lincoln Drive \\
Madison, Wisconsin 53706-1388 \\
\end{flushleft}

\appendix

\end{document}